\RequirePackage{fix-cm}
\documentclass[twocolumn]{svjour3}          % twocolumn
\smartqed  % flush right qed marks, e.g. at end of proof
\usepackage{graphicx}

\usepackage{dcolumn}
\usepackage{bm}
\usepackage[utf8]{inputenc}
\usepackage[T1]{fontenc}
\usepackage{mathptmx}
\usepackage{etoolbox}
\PassOptionsToPackage{hyphens}{url}\usepackage{hyperref}
\usepackage{amsmath}
\usepackage{amssymb}
\usepackage{amsfonts}
\usepackage{optidef}
\usepackage{algorithmic}
\usepackage{algorithm}
\usepackage{relsize}
\usepackage{appendix}

\journalname{Nonlinear Dynamics}
\begin{document}

\title{Generalized Stochastic Resilience for Early Warning Signals Based on Koopman Operator%\thanks{Grants or other notes
}


\author{Yuta Miyauchi         \and
        Masahiro Ikeda        \and 
        Yoshinobu Kawahara %etc.
}


\institute{Yuta Miyauchi \at
              Graduate School of Information Science and Technology, The University of Osaka, Osaka, Japan \\
              \email{y-miyauchi@ist.osaka-u.ac.jp}           %  \\
           \and
           Masahiro Ikeda \at
              Graduate School of Information Science and Technology, The University of Osaka, Osaka, Japan \\ RIKEN Center for Advanced Intelligence Project, Tokyo, Japan \\
              \email{ikeda@ist.osaka-u.ac.jp}
            \and
           Yoshinobu Kawahara \at
              Graduate School of Information Science and Technology, The University of Osaka, Osaka, Japan \\ RIKEN Center for Advanced Intelligence Project, Tokyo, Japan \\
              \email{kawahara@ist.osaka-u.ac.jp}
}

\date{Received: date / Accepted: date}

\maketitle

\begin{abstract}
Developing methods for detecting tipping phenomena at an early stage is an important problem in various fields such as ecology, medicine, and economics. A tipping phenomenon is characterized by a rapid transition resulting from the accumulation of small parameter changes, and is known to be related to the bifurcation theory of dynamical systems. However, few studies have examined how nonlinear properties near bifurcation points affect early warning signals (EWSs) performance. In this study, we apply the Koopman operator, which describes the time evolution of dynamical systems in an infinite-dimensional function space, to generalize stochastic resilience the theoretical basis of EWSs such as variance-based ones. As a result, we develop a novel signal capable of more accurately predicting tipping events by separately isolating stochastic fluctuations induced by noise and contributions from a continuous spectrum emerging immediately above tipping points. Our experiments indicate that our proposed method provides robust early warning detection across diverse datasets and is notably resilient to observation noise, often performing competitively with conventional indicators.
\end{abstract}

\section{Introduction}
\label{intro}
Some real-world dynamical systems undergo sudden and rapid state changes driven by variations in key parameters. These phenomena, often called \textit{tipping points} or \textit{critical transitions}\cite{scheffer2009early}, have been documented in various fields, including ecology, medicine, economics, and physics. In particular, severe environmental issues resulting from human activities, such as global warming and desertification, exhibit such tipping behaviors. Therefore, detecting \textit{early warning signals} (EWSs) is a crucial research challenge\cite{dakos2008slowing,lenton2008tipping,rockstrom2009safe}.

Among these tipping phenomena, some are directly related to the bifurcation theory in dynamical systems. Accordingly, various EWSs have been proposed based on dynamical behaviors near bifurcation points. For example, conventional EWSs include recovery rates after small perturbations\cite{van2007slow}, variance\cite{ives1995measuring}, and lag-1 autocorrelation\cite{dakos2008slowing,held2004detection}, each reflecting \textit{critical slowing down} \cite{wissel1984universal}. Critical slowing down denotes the progressively slower return to an equilibrium as a system approaches a bifurcation point. This stochastic effect appears as a decrease in resilience to noise, providing a basis for using variance as an EWS. However, most of these indicators are derived from locally linearized models in a data space, so they represent approximate linearization and do not exclusively capture the stochastic effects arising from the system’s temporal evolution.

We have therefore developed a novel approach for EWSs based on \textit{the Koopman operator}. Our method successfully generalizes stochastic resilience associated with critical slowing down. The Koopman operator, originally introduced by Koopman\cite{koopman1931hamiltonian}, is a linear operator that represents a time evolution of a system through observables. In the 2000s, it was extended to a dissipative system by Mezić\cite{mezic2005spectral}. It provides a globally linearized time evolution by lifting a time evolution of a system into an infinite-dimensional function space; it has attracted considerable attention as a powerful analytical technique for nonlinear systems\cite{brunton2022modern}. It is known that the Koopman operator can be approximated from time‐series data using Extended Dynamic Mode Decomposition (EDMD)\cite{williams2015data}. More specifically, Gaspard et al.\cite{gaspard1995spectral} show that the Liouville operator, which is the adjoint operator of the infinitesimal generator of the Koopman operator, has a continuous spectrum in the vicinity of a bifurcation point. Their research implies that learning the Koopman operator via EDMD becomes more difficult near a bifurcation point. Therefore, for predicting tipping phenomena, we focus not on how accurately the Koopman operator can be approximated, but rather on the approximation error of the Koopman operator.

In addition, our approach can be implemented via Residual Dynamic Mode Decomposition (ResDMD), an algorithm introduced by Colbrook et al.\cite{colbrook2023residual}, to quantify the estimation error of the Koopman eigenvalues and eigenfunctions. In general, ResDMD is used to detect spectral pollution arising from factors such as data quality and the accuracy of the estimation algorithm. However, it can also capture components that cannot be represented due to the projection onto a finite-dimensional space. In our method, we approximate the reconstruction error of the Koopman operator by applying ResDMD to all the estimated eigenvalue – eigenfunction pairs. Our experiments confirmed that the proposed method can indeed anticipate tipping phenomena.

Our contributions are as follows.
\begin{itemize}
    \item Nonlinear and spectral mechanisms captured in an EWS: \\We cast EWSs in the stochastic Koopman framework and show that the Koopman operator residual provides a principled signal that blows up near local bifurcation points, extending Ives’ stochastic resilience\cite{ives1995measuring} beyond local linearization (Theorem \ref{thm02}). This connects tipping prediction directly to operator spectra.
    \item Interpretability by variance-spectrum separation: \\We derive a decomposition of our EWS into a stochastic covariance term and an approximation error to the point spectrum restriction, explicitly exposing the role of the continuous spectrum that arises around bifurcation points (Eq. \eqref{eq:reskmd_decomp}). This clarifies why classical variance-based EWSs succeed or fail.
    \item Data-driven, online-computable methodology: \\Using ResDMD\cite{colbrook2023residual}, we turn the definition into a practical estimator. The average residual per eigenvalue - eigenfunction pair approximates the desired operator residual (Eq. \eqref{eq:reskmd_aprox}). We implement Window DMD with delay coordinates, enabling sliding window monitoring suitable for streaming data (Algorithms \ref{alg:ews_reskmd1}, \ref{alg:ews_reskmd2}).
    \item Practical advantages: \\The approach is pre‑training free, compatible with standard DMD toolchains, and yields diagnostic insight (how much of the dynamics cannot be represented by point spectra), which is crucial for decision making near tipping points. In the result, our methods (especially with RBF/Laplacian kernels) achieve consistently strong ROC curves and maintain performance across different parameter change rates, while a deep learning based EWS requires extensive pre‑training (Fig. \ref{fig:roc_all}).
\end{itemize}

The remainder of this paper is organized as follows. In Section 2, we review known results on EWSs for tipping phenomena induced by local bifurcations, along with the definition and properties of the stochastic Koopman operator. Based on these discussions, we introduce the stochastic residual of the Koopman mode decomposition (ResKMD) as an EWS, representing the approximation error of the Koopman operator based on its point spectra, in Section 3. Section 4 then presents a numerical algorithm based on ResDMD. In Section 5, we provide experimental examples of EWSs using both artificial and real-world datasets, comparing conventional methods with the proposed approach. Finally, Section 6 concludes the paper. Some proofs are included in the Appendix.

\section{Preliminaries}
The list of symbols used in this paper is provided in Table \ref{tab:math_symbol}.
\begin{table*}[h]
    \centering
    \caption{Table of Mathematical Symbols}
    \label{tab:math_symbol}
    \begin{tabular}{|c|c|c|c|}
        \hline \hline
        $\mathbb{N}$ & Set of natural numbers & $\mathbb{N}^{+}$ & Set of natural numbers (except for zero) \\ \hline
        $\mathbb{R}$ & Set of real numbers & $\mathbb{C}$ & Set of complex numbers \\ \hline
        $|\cdot|$ & Absolute value of $\cdot$ & $\bar{\cdot}$ & Complex conjugate of $\cdot$ \\ \hline
        $O$ & Zero matrix & $\mathrm{diag}$ & Diagonal matrix \\ \hline
        $\cdot^{T}$ & Transpose of $\cdot$ & $\cdot^*$ & Adjoint of $\cdot$ \\ \hline
        $\cdot^{-1}$ & Inverse matrix of $\cdot$ & $\cdot^{\dagger}$ & Moore–Penrose pseudoinverse matrix of $\cdot$ \\ \hline
        $\circ$ & Composition operator & $\otimes$ & Kronecker product \\ \hline
        $\mathrm{span}$ & Linear span & $L^2$ & Squared integrable function space \\ \hline
        $\|\cdot\|_2$ & Euclidean norm of $\cdot$ & $\langle \cdot_1 , \cdot_2 \rangle_2$ & Euclidean inner product between $\cdot_1$ and $\cdot_2$ \\ \hline
        $\mathcal{N}$ & Gaussian distribution & $\mathrm{Re}\cdot$ & Real part of $\cdot$ \\ \hline \hline
    \end{tabular}
\end{table*}

\subsection{Critical Slowing Down and Stochastic Resilience}
Consider a discrete-time stochastic dynamical system given by
\begin{align}
\mathbf{x}_{t+1} = \mathbf{F}(\mathbf{x}_t;\beta) + \boldsymbol{\omega}_t \triangleq \mathbf{F}_{\boldsymbol{\omega}}(\mathbf{x}_t;\beta), \label{eq:stoch_map}
\end{align}
where $t \in \mathbb{N}^+$ denotes a time step, $X\subseteq \mathbb{R}^N$ is a state space, $\mathbf{x}_t = (x_{1,t} \cdots x_{N, t})^T\in X$ is a state vector, $\Omega$ is a sample space, $\boldsymbol{\omega}_t \in \Omega$ represents a system noise, and $\beta \in \mathbb{R}$ is a bifurcation parameter. Also, $\mathbf{F}: X \rightarrow X$ is a (possibly, nonlinear) smooth vector field. We assume that the noise sequence $\{ \boldsymbol{\omega}_t \}$ is independent and identically distributed (i.i.d.) Gaussian white noise with a density $\rho(\boldsymbol{\omega}_t) = \mathcal{N}\bigl(\mathbf{0},\Sigma_{\boldsymbol{\omega}}\bigr)$, where $\rho$ is a probabilistic measure on $\Omega$. We denote $\mathrm{E}[\cdot]$ the expectation with the density $\rho$ by
\begin{align}
    \mathrm{E}[x] \triangleq \int_{\Omega} x d\rho(\boldsymbol{\omega}). \label{eq:expectation_rho}
\end{align}
Also, the variance, covariance, and autocovariance are denoted by
\begin{align}
    &\mathrm{Var}[x] \triangleq \int_{\Omega} (x - \mathrm{E}[x])^2 d\rho(\boldsymbol{\omega}), \label{eq:var_rho} \\
    &\mathrm{Cov}[x_1, x_2] \triangleq \int_{\Omega} (x_1 - \mathrm{E}[x_1]) \overline{(x_2 - \mathrm{E}[x_2])} d\rho(\boldsymbol{\omega}), \label{eq: cov_rho} \\
    & \big(\mathrm{Cov}[\mathbf{x}]\big)_{ij} \triangleq \mathrm{Cov}[x_i, x_j], \label{eq: autocov_rho}
\end{align}
where $x, x_1, x_2$ are some scalar-values, $\mathbf{x}$ is a vector-value and $x_i, x_j$ are the elements of $\mathbf{x}$.

First, by setting the noise-free limit $\Sigma_{\boldsymbol{\omega}} = O$, the system \eqref{eq:stoch_map} reduces to a purely deterministic map, which is useful to discuss the local bifurcation properties of the system. We further assume that, for $\beta < \beta^*$, where $\beta^*$ is some real value, this system admits an asymptotically stable fixed point $\mathbf{x}^* \in X$, that is, $\lim_{t\rightarrow\infty}\mathbf{x}_t=\mathbf{x}^*$. Denote by $A$ the Jacobi matrix of $\mathbf{F}$ at $\mathbf{x}^*$, i.e. 
\begin{align}
    A \triangleq \frac{\partial \mathbf{F}}{\partial \mathbf{x}} \Bigg|_{\mathbf{x} = \mathbf{x}^*}, \label{eq: matrix_jacobi}
\end{align}
whose eigenvalues $\lambda_{\mathrm{jacobi},1},\cdots,\lambda_{\mathrm{jacobi},N}$ lie strictly inside the unit circle in the complex plane:
\begin{align}
    \quad 1 > |\lambda_{\mathrm{jacobi},1}| \geq \cdots \geq |\lambda_{\mathrm{jacobi},N}| > 0. \label{eq: eigenvalue_jacobi}
\end{align}
As $\beta$ tends to $\beta^*-0$, the dominant eigenvalue $\lambda_{\mathrm{jacobi},1}$ approaches the unit circle, $|\lambda_{\mathrm{jacobi},1}| \rightarrow 1$, and $\mathbf{x}^*$ loses stability. This local stability change at $\beta = \beta^*$ is known as local bifurcation, and $\beta = \beta^*$ is called a bifurcation point\cite{kuznetsov1998elements}. Far from the bifurcation point, the trajectory of the system converges exponentially to $\mathbf{x}^*$. However, near the bifurcation point, a "resilience" toward $\mathbf{x}^*$ weakens substantially. This phenomenon is known as critical slowing down\cite{wissel1984universal,scheffer2009early}.

Next, consider a stochastic case, $\Sigma_{\boldsymbol{\omega}} \neq O$. A linearized system of \eqref{eq:stoch_map} around $\mathbf{x}^*$ yields
\begin{align}
\bar{\mathbf{x}}_{t+1} = A\bar{\mathbf{x}}_t + \boldsymbol{\epsilon}_t, \label{eq:lin_map}
\end{align}
where $\bar{\mathbf{x}}_t = \mathbf{x}_t - \mathbf{x}^*$ and $\boldsymbol{\epsilon}_t$ combine the noise $\boldsymbol{\omega}_t$ with the higher-order terms of the nonlinearity $\mathbf{F}$. Let the eigenvalue decomposition of $A$ be $A = U \Lambda U^{-1}$, where $U$ is an $N$ - dimensional invertible matrix and $\Lambda = \mathrm{diag}(\lambda_{\mathrm{jacobi},1}, \cdots, \lambda_{\mathrm{jacobi},N})$, and introduce modal coordinates,
\begin{align}
    \bar{\mathbf{z}}_t \triangleq U^{-1} \bar{\mathbf{x}}_t, \quad \mathbf{e}_t \triangleq U^{-1} \boldsymbol{\epsilon}_t. \label{eq:coordinate_modal}
\end{align}
Focusing on the dominant mode, $\bar{z}_{1, t}$ satisfies
\begin{align}
\bar{z}_{1,t+1} = \lambda_{\mathrm{jacobi},1}\bar{z}_{1,t} + e_{1,t}. \label{eq:dom_mode}
\end{align}
Furthermore, we assume that Eq. \eqref{eq:dom_mode} is a stationary process in $\bar{z}_{1, t}$, and that $\bar{z}_{1, t}$ and $e_{1, t}$ are orthogonal. Under these assumptions for $\bar{z}_{1, t}$, its variance is given by
\begin{align}
\mathrm{Var}[\bar{z}_{1,t}]
= \frac{\mathrm{Var}[e_{1,t}]}{1 - |\lambda_{\mathrm{jacobi},1}|^2}. \label{eq:var_diverge}
\end{align}
Hence, as $|\lambda_{\mathrm{jacobi},1}| \rightarrow 1$, $\mathrm{Var}[\bar{z}_{1, t}]$ diverges to infinity\cite{ives1995measuring}, reflecting the weakened "resilience" via critical slowing down.

\subsection{Stochastic Koopman Operator}
Let $g: \mathbb{R}^N \rightarrow \mathbb{C}$ be an observable. Its time evolution satisfies
\begin{align}
g(\mathbf{x}_{t+1}) &= g\bigl(\mathbf{F}_{\boldsymbol{\omega}}(\mathbf{x}_t;\beta)\bigr) = (g \circ \mathbf{F}_{\boldsymbol{\omega}})(\mathbf{x}_t;\beta) \label{eq:koop_map}.
\end{align}
Since $g \circ \mathbf{F}_{\boldsymbol{\omega}}$ is a random variable, its expectation is called the stochastic Koopman operator\cite{vcrnjaric2020koopman},
\begin{align}
U_{(1)}g(\mathbf{x}) &\triangleq \mathrm{E}\bigl[g\circ\mathbf{F}_{\boldsymbol{\omega}}(\mathbf{x};\beta)\bigr] \label{eq:koop_sdef}.
\end{align}
We consider $U_{(1)}$ the operator on the squared integrable function space $L^2(X, \mu)$, where $\mu$ is a positive measure on $X$. Here we define the norm $\|\cdot\|_{L^2(X,\mu)}$ and the inner product $\langle \cdot,\cdot\rangle_{L^2(X,\mu)}$ on $L^2(X,\mu)$ as
\begin{align}
    \| g \|_{L^2(X,\mu)}^2 &\triangleq \int_X |g(\mathbf{x})|^2 d\mu(\mathbf{x}), \label{eq:norm_obs} \\ 
    \langle g_1,g_2\rangle_{L^2(X,\mu)}&\triangleq\int_{X}g_1(\mathbf{x})\overline{g_2(\mathbf{x})}d\mu(\mathbf{x}), \label{eq:inner_obs}
\end{align}
where $g, g_1, g_2 \in L^2(X, \mu)$. In $\Sigma_{\boldsymbol{\omega}} = O$, one recovers the deterministic Koopman operator:
\begin{align}
U_{(1)}g = g \circ \mathbf{F}. \label{eq:koop_ddef}
\end{align}

Since $U_{(1)}$ is a linear operator, its eigenvalue and eigenfunction can be defined, and their various estimation methods have been extensively studied\cite{mauroy2020koopman}. Let $\lambda \in \mathbb{C}$ be a Koopman eigenvalue of $U_{(1)}$ with the corresponding Koopman eigenfunction $\phi_{\lambda} \in L^2(X, \mu)$. Then, we have
\begin{align}
U_{(1)}\phi_\lambda = \lambda\phi_\lambda. \label{eq:koop_eig}
\end{align}
If the system \eqref{eq:stoch_map} possesses a hyperbolic fixed point (not in the vicinity of the unit circle), the Poincaré linearization theorem guaranties a nonlinear coordinate that linearizes the flow near the fixed point. Hence, in a suitable function space, $U_{(1)}$ has only countably infinite point spectra $\{\lambda_1,\lambda_2,\cdots\}$ with $|\lambda_1|\geq|\lambda_2|\geq\cdots$\cite{mauroy2016global}. In particular, in $\Sigma_{\boldsymbol{\omega}} = O$, the flow of the system can be expressed as
\begin{align}
U_{(1)}^t g(\mathbf{x})
= \sum_{n=1}^\infty \lambda_n^tv_{\lambda_n}\phi_{\lambda_n}(\mathbf{x}), \label{eq:koop_decomp}
\end{align}
where $v_{\lambda_n} \in \mathbb{C}$ is the Koopman mode associated with $\lambda_n$. In contrast, when the fixed point is nonhyperbolic, namely some Koopman eigenvalues are dense on the unit circle, $U_{(1)}$ exhibits a continuous spectrum\cite{gaspard1995spectral,mauroy2016global}. This emergence of the continuous spectrum may be a hallmark of local bifurcation.

Moreover, since the action of the stochastic Koopman operator  is the expectation of the random variable $g \circ \mathbf{F}_{\boldsymbol{\omega}}$, it is natural to define its variance by analogy with broader stochastic frameworks. To define the variance, we introduce a batched observable $h: \mathbb{R}^N \times \mathbb{R}^N \rightarrow \mathbb{C}$ and the batched Koopman operator as follows,
\begin{align}
U_{(2)}h \triangleq \mathrm{E}\bigl[h\bigl(\mathbf{F}_{\boldsymbol\omega},\mathbf{F}_{\boldsymbol\omega'}\bigr)\bigr],\label{eq:koop_batched}
\end{align}
which is firstly introduced by Colbrook et al\cite{colbrook2024beyond}. Then, the variance along the trajectory through $\mathbf{x}$ is yielded by
\begin{align}
\mathrm{Var}\bigl[g\circ\mathbf{F}_{\boldsymbol\omega}(\mathbf{x};\beta)\bigr]
&= U_{(2)}\bigl(g\otimes\overline{g}\bigr)(\mathbf{x},\mathbf{x})
- \bigl\lvert U_{(1)}g(\mathbf{x})\bigr\rvert^2. \label{eq:kooop_localvar}
\end{align}
Integrating (\ref{eq:kooop_localvar}) over the state space, the global variance is defined as follows
\begin{align}
\mathrm{Var}[g\circ\mathbf{F}_{\boldsymbol\omega}]
&\triangleq \int_X \mathrm{Var}\bigl[g\circ\mathbf{F}_{\boldsymbol\omega}(\mathbf{x};\beta)\bigr]d\mu(\mathbf{x}). \label{eq:koop_globalvar}
\end{align}
This variance \eqref{eq:koop_globalvar} satisfies the following theorem\cite{colbrook2024beyond}, which provides a foundation for estimating the variance given by \eqref{eq:koop_globalvar}. 
\begin{theorem} 
Let $g_1, g_2 \in L^2(X, \mu)$. Then, 
\begin{align}
&\mathrm{E}\bigl[\| g_1\circ\mathbf{F}_{\boldsymbol\omega}+g_2\|_{L^2(X, \mu)}^2\bigr] \notag \\ 
&\qquad= \| U_{(1)}g_1 + g_2\|_{L^2(X, \mu)}^2 + \mathrm{Var}[g_1\circ\mathbf{F}_{\boldsymbol\omega}]. \label{eq:single_biasvariance}
\end{align}
\label{thm01} 
\end{theorem}
If we use $g_2$ as an estimate of the eigenfunction of the Koopman operator, the left-hand side of \eqref{eq:single_biasvariance} represents the integrated mean squared error (IMSE).

Also, covariance admits the similar formulation. For $g_1, g_2 \in L^2(X, \mu)$, the covariance operator is defined by
\begin{align}
&\mathrm{Cov}[g_1\circ\mathbf{F}_{\boldsymbol\omega},g_2\circ\mathbf{F}_{\boldsymbol\omega}] \notag \\ 
& \qquad \triangleq \int_X \mathrm{Cov}\big[g_1\circ\mathbf{F}_{\boldsymbol\omega}(\mathbf{x};\beta), g_2\circ\mathbf{F}_{\boldsymbol\omega}(\mathbf{x};\beta)\big] d\mu(\mathbf{x}). \label{eq:Koop_covariance}
\end{align}

\section{Koopman Framework for Early Warning Signals}

\begin{figure}[t]
    \centering
    \includegraphics[width=0.5\textwidth]{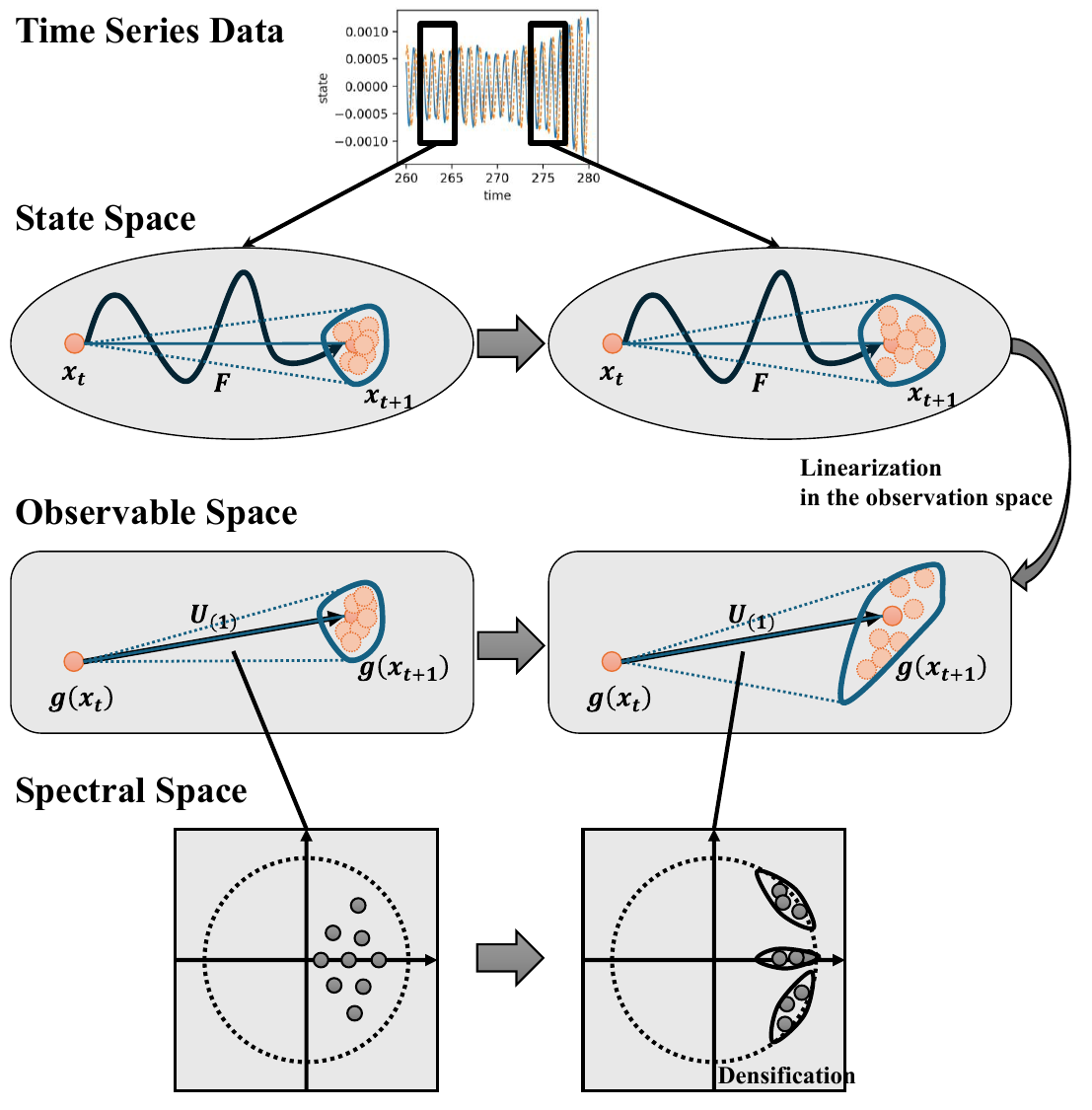}
    \caption{Conceptual diagram of ResKMD}
    \label{fig:ex_reskmd}
\end{figure}

In this section, we establish a novel Koopman framework for adopting an EWS. Especially, we focus on an approximator of the stochastic Koopman operator $U_{(1)}$ by numerical estimation methods such as Extended Dynamic Mode Decomposition (EDMD)\cite{williams2015data}. Then we consider a restriction of the Koopman operator in the following invariant subspace in $L^2(X,\mu)$, 
\begin{align}
    \Phi_m = \mathrm{span}\{ \phi_{\lambda_1}, \cdots, \phi_{\lambda_m} \}, \label{eq:invariant_subspace}
\end{align}
where $m\in \mathbb{N}^+$ and $\phi_{\lambda_n} \in L^2(X, \mu)$ for $n = 1, \cdots, m$ are the Koopman eigenfunctions with the Koopman eigenvalues $\lambda_n$. The restricted stochastic Koopman operator is denoted by $U_{(1), \Phi_m}$. Then, the time evolution of an $M$-dimensional vector-valued observable $\mathbf{g}: X \rightarrow \mathbb{C}^M$ can be expressed by\cite{takeishi2017subspace}
\begin{align}
    U_{(1), \Phi_m}\mathbf{g}(\mathbf{x}) = \sum_{n=1}^m \lambda_n \mathbf{v}_{\lambda_n}\phi_{\lambda_n}(\mathbf{x}). \label{eq:restricted_koop}
\end{align}
Here, a vector-valued observable is $M$ collections of scalar-valued observables in $L^2(X, \mu)$. This expression \eqref{eq:restricted_koop} approximates the true Koopman operator through its dominant point spectra. However, for the emergence of a continuous spectrum at the bifurcation point, the corresponding part of $U_{(1)}$ cannot be captured by $U_{(1), \Phi_m}$ and this phenomenon can be considered as an EWS. To quantify this missing part, we define the stochastic residual of the Koopman mode decomposition (ResKMD) by
\begin{align}
    \operatorname{res}\bigl[U_{(1)},U_{(1),\Phi_m}\bigr]^2 \triangleq \mathrm{E}\Bigl[\bigl\|\mathbf{g}\circ\mathbf{F}_{\boldsymbol{\omega}}-U_{(1),\Phi_m}\mathbf{g}\bigr\|_{L^2(X, \mu)^M}^2\Bigr], \label{eq:def_reskmd}
\end{align}
where $\| \cdot\|_{L^2(X, \mu)^M}$ is the norm on $L^2(X, \mu)^M$ defined by
\begin{align}
    \| \mathbf{g}\|_{L^2(X, \mu)^M} \triangleq \Bigg(\sum_{n=1}^M \| g_n \|_{L^2(X, \mu)}^2 \Bigg)^{\frac{1}{2}} \label{eq:norm_vectorobs}
\end{align}
for $\mathbf{g}=(g_1 \cdots g_M)^T\in L^2(X, \mu)^M$.
We examine properties of ResKMD near the bifurcation point. To decompose the residual, we extend Theorem \ref{thm01} for vector-valued observables.
\begin{corollary}
    Let $\mathbf{g}_1, \mathbf{g}_2: X\rightarrow\mathbb{C}^M$ be vector-valued observables. Then,
    \begin{align}
    &\mathrm{E}\bigl[\|\mathbf{g}_1\circ\mathbf{F}_{\boldsymbol{\omega}} + \mathbf{g}_2\|_{L^2(X, \mu)^M}^2\bigr] \notag \\
    & \qquad = \|U_{(1)}\mathbf{g}_1 + \mathbf{g}_2\|_{L^2(X, \mu)^M}^2 + \operatorname{Tr}\bigl(\mathrm{Cov}[\mathbf{g}_1\circ\mathbf{F}_{\boldsymbol\omega}]\bigr), \label{eq:multi_biasvariance}
    \end{align}
    where the autocovariance operator of a vector-valued observable $\mathbf{g} = (g_1, \cdots, g_M)^T$ is defined by,
    \begin{align}
        \bigl(\mathrm{Cov}[\mathbf{g}\circ\mathbf{F}_{\boldsymbol\omega}]\bigr)_{ij} \triangleq \mathrm{Cov}\bigl[g_i\circ\mathbf{F}_{\boldsymbol{\omega}}, g_j\circ\mathbf{F}_{\boldsymbol{\omega}}\bigr], \quad 1 \leq i,j \leq M, \notag
    \end{align}
    and $\operatorname{Tr}$ is a trace operator.
    Especially, choosing $\mathbf{g}_1 = \mathbf{g}$ and $\mathbf{g}_2 = -U_{(1), \Phi_m}\mathbf{g}$ yields the decomposition of ResKMD below,
\begin{align}
&\operatorname{res}\bigl[U_{(1)},U_{(1),\Phi_m}\bigr]^2 \notag \\
& \qquad = \|U_{(1)}\mathbf{g} - U_{(1),\Phi_m}\mathbf{g}\|_{L^2(X, \mu)^M}^2 + \operatorname{Tr}\bigl(\mathrm{Cov}[\mathbf{g}\circ\mathbf{F}_{\boldsymbol\omega}]\bigr). \label{eq:reskmd_decomp}
\end{align}
    \label{cor01} 
\end{corollary}

Let us focus on the first term on the right-hand side of Eq. \eqref{eq:reskmd_decomp}. This quantity means the squared error between the true Koopman operator $U_{(1)}$ and its restriction $U_{(1), \Phi_m}$. For the system \eqref{eq:stoch_map} with a hyperbolic and asymptotically stable fixed point, the following proposition holds.
\begin{proposition}
    Let $\mathbf{g}: X\rightarrow\mathbb{C}^M$ be a vector-valued smooth observable, and assume further that the set of Koopman eigenfunctions forms an orthonormal basis. Then, 
    \begin{align}
        \lvert \lambda_{m+1} \rvert^2 \| \mathbf{v}_{\lambda_{m+1}} \|_2^2 &\le \bigl\|U_{(1)}\mathbf{g}-U_{(1),\Phi_m}\mathbf{g}\bigr\|_{L^2(X, \mu)^M}^2 \notag \\
        & \qquad\le \lvert \lambda_{m+1} \rvert^2 \|\mathbf{g}\|_{L^2(X, \mu)^M}^2. \label{eq:hypb_resbounds}
    \end{align}
    \label{prop01} 
\end{proposition}
Sufficiently far from the bifurcation point, the norm of the Koopman eigenfunctions and modes vary little. Accordingly, as the $(m + 1)$-th largest eigenvalue gradually approaches the unit cycle, Eq. \eqref{eq:hypb_resbounds} implies that the squared error between $U_{(1)}$ and $U_{(1), \Phi_m}$ increases.

Next, we focus on the second term on the right-hand side of Eq. \eqref{eq:reskmd_decomp}, which captures the stochastic influence in the evolution of observables. Since the autocovariance operator is a Hilbert-Schmidt operator, we obtain
\begin{align}
    \operatorname{Tr}\bigl(\mathrm{Cov}[\mathbf{g}\circ\mathbf{F}_{\boldsymbol\omega}]\bigr) = \sum_{n=1}^\infty\lambda_{\mathrm{cov},n}, \notag
\end{align}
where $\{\lambda_{\mathrm{cov},n}\}_{n=1}^\infty$ are the eigenvalues of $\operatorname{Cov[g\circ\mathbf{F}_{\boldsymbol\omega}]}$. For each $n\in \mathbb{N}$, the eigenvalue $\lambda_{\mathrm{cov},n}$ quantifies the variance injected by the noise along the Koopman eigenfunction $\phi_{\lambda_n}$. For every such $\phi_{\lambda_n}$, we have the global evolution law,
\begin{align}
    \phi_{\lambda_n}(\mathbf{x}_{t+1}) = \lambda_n \phi_{\lambda_n}(\mathbf{x}_t) + \eta_{n, t}, \label{eq:eig_evolution}
\end{align}
where the random variable $\eta_{n, t} = \phi_{\lambda_n}\circ\mathbf{F}_{\boldsymbol{\omega}}(\mathbf{x};\beta) - \lambda_n\phi_{\lambda_n}(\mathbf{x})$ represents the influence of noise $\boldsymbol{\omega}_t$ and the variance of $\phi_{\lambda_n}$ is given by $\lambda_{\mathrm{cov},n}$. We assume that $\eta_{n, t}$ is independent of $\mathbf{x}_t$, repeating the argument that leads to Eq. \eqref{eq:var_diverge} yields the following result.
\begin{align}
    \operatorname{Tr}\bigl(\mathrm{Cov}[\mathbf{g}\circ\mathbf{F}_{\boldsymbol\omega}]\bigr) = \sum_{n=1}^\infty \frac{\mathrm{Var}[\eta_{n,t}]}{1 - |\lambda_n|^2}. \label{eq:var_generalize} 
\end{align}
Hence, as the system approaches the bifurcation point and the absolute values of the dominant eigenvalues $\lambda_n$ tend to $1$, the second term in \eqref{eq:reskmd_decomp} increases.

In summary, the behavior of ResKMD in a neighborhood of the bifurcation point can be formulated as the following theorem.
\begin{theorem} 
Consider the stochastic dynamical system \eqref{eq:stoch_map}. If the system has a hyperbolic and asymptotically stable fixed point, then
\begin{align}
    \operatorname{res}\bigl[U_{(1)},U_{(1),\Phi_m}\bigr]^2 \rightarrow \infty \quad \mathrm{as} \quad \beta \rightarrow \beta^*. \label{eq:ews_reskmd} 
\end{align}
\label{thm02} 
\end{theorem}
Theorem \ref{thm02} can be regarded as an extension of Ives' approach to the resilience for the EWS\cite{ives1995measuring}. Ives' approach focused on the local behavior of the linearized system near the bifurcation point. However, our Koopman-based approach can be considered nonlinear features that cannot be embedded in the linear term. (The graphical explanation is shown by Fig. \ref{fig:ex_reskmd}.)

\section{Numerical Method}
\begin{algorithm*}[t]
\caption{Early Warning Signal by Calculating ResKMD with Exact DMD}
\label{alg:ews_reskmd1}
\begin{algorithmic}[1]
\REQUIRE Snapshot data $\{ \mathbf{x}^{(t)}, \mathbf{y}^{(t)} \}_{t=1}$, quadrature weights $W$, positive integer $T_{\mathrm{window}}, d_{\mathrm{hankel}}, r$.
\FOR{$i = 1, 2, \cdots$}
\STATE Extract window data $\{ \mathbf{x}^{(t)}, \mathbf{y}^{(t)} \}_{t=i}^{i + T_{\mathrm{window}}}$ and extend to $\{ \mathbf{x}_{\mathrm{hankel}}^{(t)}, \mathbf{y}_{\mathrm{hankel}}^{(t)} \}_{t=i}^{i + T_{\mathrm{window}} - d_{\mathrm{hankel}} + 1}$ with delay coordinate ($\mathbb{R}^N \rightarrow \mathbb{R}^{d_{\mathrm{hankel}}N}$).
\STATE Set $\Psi_{\mathrm{DMD}}(\mathbf{x}_{\mathrm{hankel}}) = [\mathbf{e}_1^* \mathbf{x}_{\mathrm{hankel}} \cdots \mathbf{e}_{d_{\mathrm{hankel}} \times N}^* \mathbf{x}_{\mathrm{hankel}}]$, where $\mathbf{e}_i$ is the $i$th unit vector.
\STATE Compute a truncated SVD
\begin{align}
\frac{1}{T_{\mathrm{window}} - d_{\mathrm{hankel}} + 1}(\Psi_{\mathrm{DMD}}(\mathbf{x}_{\mathrm{hankel}}^{(i)})^T \cdots \Psi_{\mathrm{DMD}}(\mathbf{x}_{\mathrm{hankel}}^{(i + T_{\mathrm{window}} - d_{\mathrm{hankel}} + 1)})^T)^T \approx U_r \Sigma_r V_r^*. \notag
\end{align}
\STATE Set two matrices
\begin{align}
\Psi_X = \begin{pmatrix}
        \Psi_{\mathrm{DMD}}(\mathbf{x}_{\mathrm{hankel}}^{(i)}) \\
        \vdots \\
        \Psi_{\mathrm{DMD}}(\mathbf{x}_{\mathrm{hankel}}^{(i + T_{\mathrm{window}} - d_{\mathrm{hankel}} + 1)})
    \end{pmatrix} V_r \Sigma_r^{\dagger}, \quad 
\Psi_Y = \begin{pmatrix}
        \Psi_{\mathrm{DMD}}(\mathbf{y}_{\mathrm{hankel}}^{(i)}) \\
        \vdots \\
        \Psi_{\mathrm{DMD}}(\mathbf{y}_{\mathrm{hankel}}^{(i + T_{\mathrm{window}} - d_{\mathrm{hankel}} + 1)})
    \end{pmatrix} V_r \Sigma_r^{\dagger}. \notag
\end{align}
\STATE Solve eigendecomposition
\begin{align}
(\Psi_X^*W\Psi_X)^{\dagger}(\Psi_X^*W\Psi_Y) \boldsymbol{\xi}_{\lambda_n} = \lambda_n \boldsymbol{\xi}_{\lambda_n}, \quad (n = 1, \cdots, r). \notag
\end{align}
\STATE Compute $\mathrm{res}[\lambda_n, \phi_{\lambda_n}]$ by Eq.\eqref{eq:resdmd_bydata} for all $n$.
\STATE ResKMD is given by \eqref{eq:reskmd_aprox}.
\ENDFOR
\RETURN ResKMD calculated at each time window.
\end{algorithmic}
\end{algorithm*}
\begin{algorithm*}[t]
\caption{Early Warning Signal by Calculating ResKMD with Kernel EDMD}
\label{alg:ews_reskmd2}
\begin{algorithmic}[1]
\REQUIRE Snapshot data $\{ \mathbf{x}^{(t)}, \mathbf{y}^{(t)} \}_{t=1}$, quadrature weights $W$, positive-definite kernel function $S: X \times X \rightarrow \mathbb{R}$ and positive integer $T_{\mathrm{window}}, d_{\mathrm{hankel}}, r$.
\FOR{$i = 1, 2, \cdots$}
\STATE Extract window data $\{ \mathbf{x}^{(t)}, \mathbf{y}^{(t)} \}_{t=i}^{i + T_{\mathrm{window}}}$ and extend to $\{ \mathbf{x}_{\mathrm{hankel}}^{(t)}, \mathbf{y}_{\mathrm{hankel}}^{(t)} \}_{t=i}^{i + T_{\mathrm{window}} - d_{\mathrm{hankel}} + 1}$ with delay coordinate ($\mathbb{R}^N \rightarrow \mathbb{R}^{d_{\mathrm{hankel}}N}$).
\STATE Generate gram matrices $\Psi_X, \Psi_Y$ for $\{ \mathbf{x}_{\mathrm{hankel}}^{(t)}, \mathbf{y}_{\mathrm{hankel}}^{(t)} \}_{t=i}^{i + T_{\mathrm{window}} - d_{\mathrm{hankel}} + 1}$ with kernel $S$ and compute $r$-rank approximated SVD of $\Psi_X$ and $\tilde{\mathbb{K}}$ as follows,
\begin{align}
    \sqrt{W} \Psi_X \Psi_X^* \sqrt{W} = U_r \Sigma_r^2 U_r^*, \quad\tilde{\mathbb{K}} = (\Sigma_r^{\dagger}U_r^*)(\sqrt{W} \Psi_Y \Psi_X^* \sqrt{W})(U_r\Sigma_r^{\dagger}). \notag
\end{align}
\STATE Compute the eigenvalues of $\tilde{\mathbb{K}}$ and stack the corresponding eigenvectors column-by-column into $Z \in \mathbb{C}^{(T_{\mathrm{window}} - d_{\mathrm{hankel}} + 1) \times r}$.
\STATE Apply a QR decomposition to orthogonalize $Z$ with respect to $Q = [Q_1, \cdots, Q_r]$.
\STATE Set the dictionary as follows,
\begin{align}
    \psi_{\mathrm{selected},j}(\mathbf{x}_{\mathrm{hankel}}) = \Big[S\big(\mathbf{x}_{\mathrm{hankel}}, \mathbf{x}_{\mathrm{hankel}}^{(i)}\big), \cdots, S\big(\mathbf{x}_{\mathrm{hankel}}, \mathbf{x}_{\mathrm{hankel}}^{(i + T_{\mathrm{window}} - d_{\mathrm{hankel}} + 1)}\big)\Big](U_r \Sigma_r^{\dagger})Q_j, \quad 1 \leq j \leq r. \notag
\end{align}
\STATE Solve eigenvalue decomposition
\begin{align}
(\Psi_{\mathrm{selected},X}^*W\Psi_{\mathrm{selected},X})^{\dagger}(\Psi_{\mathrm{selected},X}^*W\Psi_{\mathrm{selected},Y}) \boldsymbol{\xi}_{\lambda_n} = \lambda_n \boldsymbol{\xi}_{\lambda_n}, \quad (n = 1, \cdots, r). \notag
\end{align}
\STATE Compute $\mathrm{res}[\lambda_n, \phi_{\lambda_n}]$ by Eq.\eqref{eq:resdmd_bydata} for all $n$.
\STATE ResKMD is given by \eqref{eq:reskmd_aprox}.
\ENDFOR
\RETURN ResKMD calculated at each time window.
\end{algorithmic}
\end{algorithm*}

\subsection{Extended Dynamic Mode Decomposition}
Suppose that we are given a snapshot dataset $\{\mathbf{x}^{(i)}, \mathbf{y}^{(i)}\}_{i=1}^T$, where $T \in \mathbb{N}$ is the number of snapshot data and each snapshot pair $\{ \mathbf{x}^{(i)},\mathbf{y}^{(i)}\}$ satisfies
\begin{align}
    \mathbf{y}^{(i)} = \mathbf{F}_{\boldsymbol{\omega}}(\mathbf{x}^{(i)}; \beta). \label{eq:eq_snapshot}
\end{align}

In EDMD\cite{williams2015data}, we need to prepare a dictionary of nonlinear functions $\Psi(\mathbf{x}) = (\psi_1(\mathbf{x}) \cdots \psi_M(\mathbf{x})) \in \mathbb{C}^{1 \times M}$ for $\mathbf{x}\in X$, and any new observable $g \in L^2(X, \mu)$ is expressed via a coefficient vector $\boldsymbol{\zeta} = (\zeta_1 \cdots \zeta_M)^T \in \mathbb{C}^M$ as
\begin{align}
    g(\mathbf{x}) = \Psi(\mathbf{x})\boldsymbol{\zeta}. \label{eq:observable_by_dictionary}
\end{align}
Under this representation, the Koopman operator $U_{(1)}$ is represented by
\begin{align}
    U_{(1)}g(\mathbf{x})&=\Psi\bigl(\mathbf{F}_{\boldsymbol{\omega}}(\mathbf{x}; \beta)\bigr)\boldsymbol{\zeta} = \Psi(\mathbf{x})\mathbb{K}\boldsymbol{\zeta} + \mathbf{R}(\boldsymbol{\zeta}, \mathbf{x}), \label{eq:koopman_approx}
\end{align}
where
\begin{align}
    \mathbf{R}(\boldsymbol{\zeta}, \mathbf{x})\triangleq \Psi\bigl(\mathbf{F}_{\boldsymbol{\omega}}(\mathbf{x}; \beta)\bigr)\boldsymbol{\zeta}-\Psi(\mathbf{x})\mathbb{K}\boldsymbol{\zeta}. \label{eq:error_kooparox}
\end{align}
Minimizing Eq. \eqref{eq:error_kooparox} with respect to $\mathbb{K}$ yields the EDMD-based approximation of the Koopman operator\cite{williams2015data}, namely:
\begin{align}
    \mathbb{K} \triangleq (\Psi_X^*W\Psi_X)^{\dagger}(\Psi_X^*W\Psi_Y), \label{eq:dmd_matrix}
\end{align}
where $Z^\dagger$ is the Moore-Penrose pseudoinverse of the matrix $Z$, $W = \mathrm{diag}(w_1, \cdots, w_T)$ is a weight matrix that quantifies the relative importance assigned to each time step and
\begin{align}
    &\Psi_X \triangleq \Big(\Psi(\mathbf{x}^{(1)})^T \cdots \Psi(\mathbf{x}^{(T)})^T \Big)^T \in \mathbb{C}^{T \times M}, \label{eq:datamat_X} \\
    &\Psi_Y \triangleq \Big(\Psi(\mathbf{y}^{(1)})^T \cdots \Psi(\mathbf{y}^{(T)})^T \Big)^T \in \mathbb{C}^{T \times M}. \label{eq:datamat_XY}
\end{align}

As $T \rightarrow \infty$, we have
\begin{align}
    \lim_{T\rightarrow\infty}(\Psi_X^*W\Psi_X)_{ij} &= \langle \psi_j, \psi_i \rangle_{L^2(X, \mu)} \label{eq:galerkin_X} \\
    \lim_{T\rightarrow\infty}(\Psi_X^*W\Psi_Y)_{ij} &= \langle U_{(1)}\psi_j, \psi_i \rangle_{L^2(X, \mu)}. \label{eq:galerkin_XY}
\end{align}
Therefore, as T tends to $\infty$, $\mathbb{K}$ becomes the true Koopman operator $U_{(1)}$.

\subsection{ResKMD by Residual Dynamic Mode Decomposition}
Colbrook et al.\cite{colbrook2023residual} proposed an evaluation framework called Residual Dynamic Mode Decomposition (ResDMD), which evaluates the reliability of the candidate of eigenvalue – eigenfunction pairs by computing $U_{(1)}^*U_{(1)}$. For any pair $\{ \lambda, \phi_\lambda \}$, ResDMD computes the following residual, denoted by $\mathrm{res}[\lambda, \phi_{\lambda}]^2$:
\begin{align}
    \mathrm{res}[\lambda, \phi_{\lambda}]^2 \triangleq \frac{\|\phi_{\lambda} \circ \mathbf{F}_{\boldsymbol{\omega}} - \lambda \phi_{\lambda}\|_{L^2(X, \mu)}^2}{\|\phi_{\lambda}\|_{L^2(X, \mu)}^2}. \label{eq:def_resdmd}
\end{align}

Next, similarly to Eqs. \eqref{eq:galerkin_X} and \eqref{eq:galerkin_XY}, we have
\begin{align}
    \lim_{T\rightarrow\infty}(\Psi_Y^*W\Psi_Y)_{ij} = \langle U_{(1)}\psi_j, U_{(1)}\psi_i \rangle_{L^2(X, \mu)}. \label{eq:galerkin_Y}
\end{align}
Then, letting $\phi_{\lambda}(\mathbf{x}) = \Psi(\mathbf{x})\boldsymbol{\xi}_{\lambda}$, where $\boldsymbol{\xi}_{\lambda}$ is the eigenvector of the matrix $\mathbb{K}$ corresponding to the eigenvalue $\lambda$, the residual \eqref{eq:def_resdmd} is approximated by
\begin{align}
    &{\small \mathrm{res}[\lambda, \phi_{\lambda}]^2} \notag \\ 
    &\qquad {\small \approx \frac{\boldsymbol{\xi}_{\lambda}^*[\Psi_Y^*W\Psi_Y-\lambda(\Psi_X^*W\Psi_Y)^* - \bar{\lambda}(\Psi_X^*W\Psi_Y) + |\lambda|^2\Psi_X^*W\Psi_X]\boldsymbol{\xi}_{\lambda}}{\boldsymbol{\xi}_{\lambda}^*[\Psi_X^*W\Psi_X]\boldsymbol{\xi}_{\lambda}}. \label{eq:resdmd_bydata}}
\end{align}
In the stochastic setting, the quantity \eqref{eq:def_resdmd} computed by ResDMD is equal to Eq. \eqref{eq:single_biasvariance} with $g_1 = \phi_{\lambda}$ and $g_2 = - \lambda \phi_{\lambda}$. This fact plays an important role in the numerical approximation of ResKMD. In fact, for candidate pairs $\{ \lambda_n, \phi_{\lambda_n} \}_{n = 1}^m$, ResKMD can be approximated by the following expression using ResDMD (also see the Appendix):
\begin{align}
    \mathrm{res}[U_{(1)}, U_{(1), \Phi_m}]^2 \approx \frac{1}{m} \sum_{n=1}^m \mathrm{res}[\lambda_n, \phi_{\lambda_n}]^2.  \label{eq:reskmd_aprox}
\end{align}

\subsection{Online Estimation of ResKMD}
In order to compute ReKMD for an EWS, we need to calculate DMD in an online manner. These are two main approaches for this purpose; Weighted Dynamic Mode Decomposition (Weighted DMD), which incorporates a forgetting coefficient for past data, and Windowed Dynamic Mode Decomposition (Windowed DMD), which applies DMD at fixed time intervals (time window)\cite{zhang2019online,alfatlawi2019incremental}. In our experiments, we choose an online approach using Window DMD because of computational complexity. Over a rolling window, we can track the temporal evolution of ResKMD. In addition, we adopt time-delayed coordinates applied to observed data to enhance the learning accuracy of Hankel DMD\cite{takens2006detecting}. 

Finally, Algorithm \ref{alg:ews_reskmd1} and Algorithm \ref{alg:ews_reskmd2} above present the algorithms to calculate the EWS based on ResKMD. Algorithm \ref{alg:ews_reskmd1} uses ResDMD by estimating the Koopman eigenvalue - eigenfunction pairs with exact DMD. Algorithm \ref{alg:ews_reskmd2} estimates the pairs with kernel EDMD. Note that these algorithms calculate the residual \eqref{eq:resdmd_bydata} following the procedure proposed by Colbrook et al.\cite{colbrook2023residual}.

\section{Experiment}

\begin{figure*}[t]
    \centering
    \includegraphics[width=\textwidth]{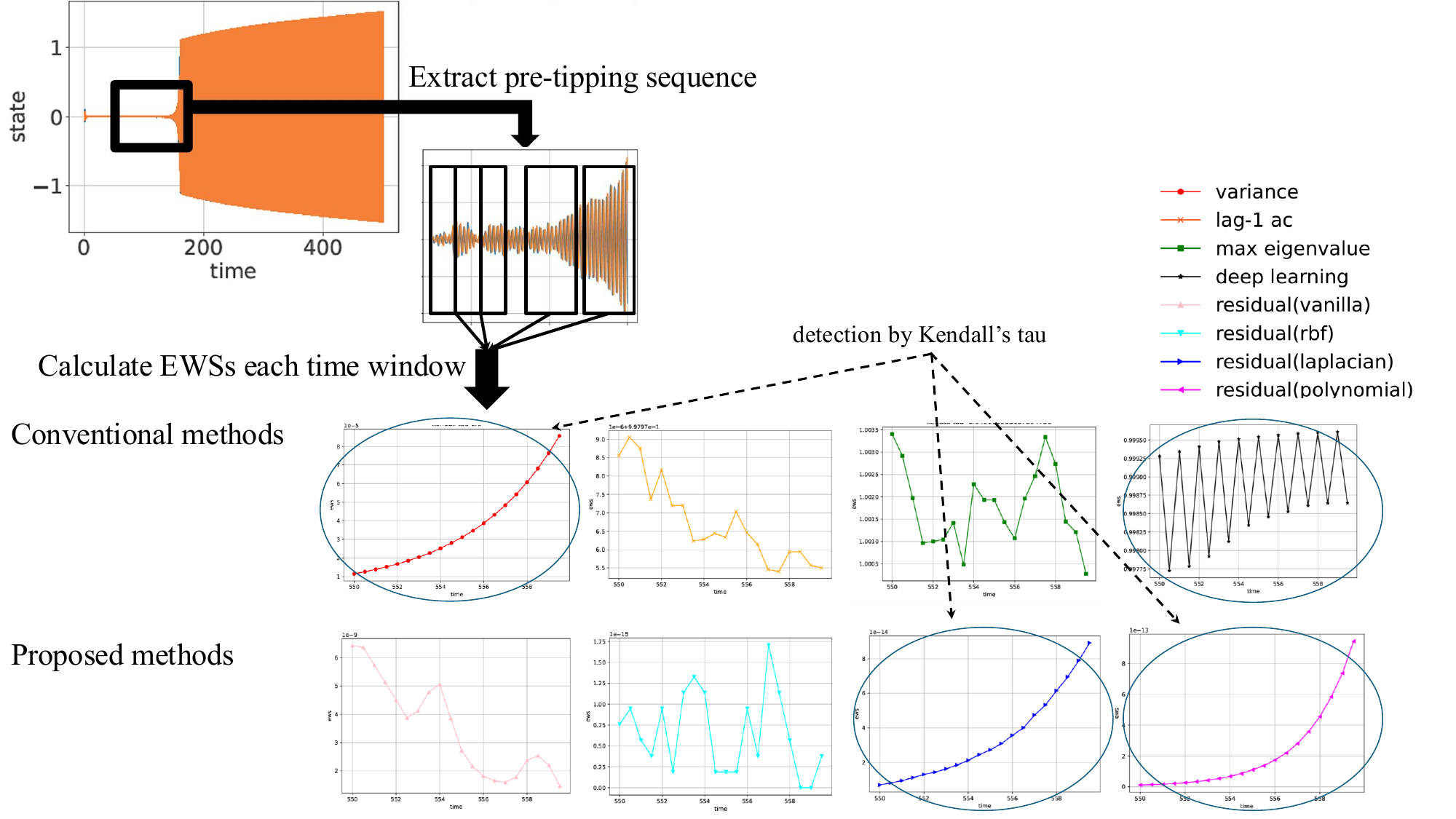}
    \caption{Example of Subcritical Hopf Bifurcation: EWSs are computed by sliding a time window along the data. Moreover, suitable EWSs are expected to increase as the system approaches the bifurcation point. Namely, this implies that Kendall’s $\tau$ between the EWSs and time should approach $1$. (Residual(vanilla) means ResKMD calculated by Exact DMD. Also, residual(rbf, laplacian, polynomial) are ResKMD calculated by EDMD with RBF, Laplacian, and polynomial kernel.)}
    \label{fig:ex_hopf}
\end{figure*}

\begin{figure*}[t]
    \centering
    \includegraphics[width=\textwidth]{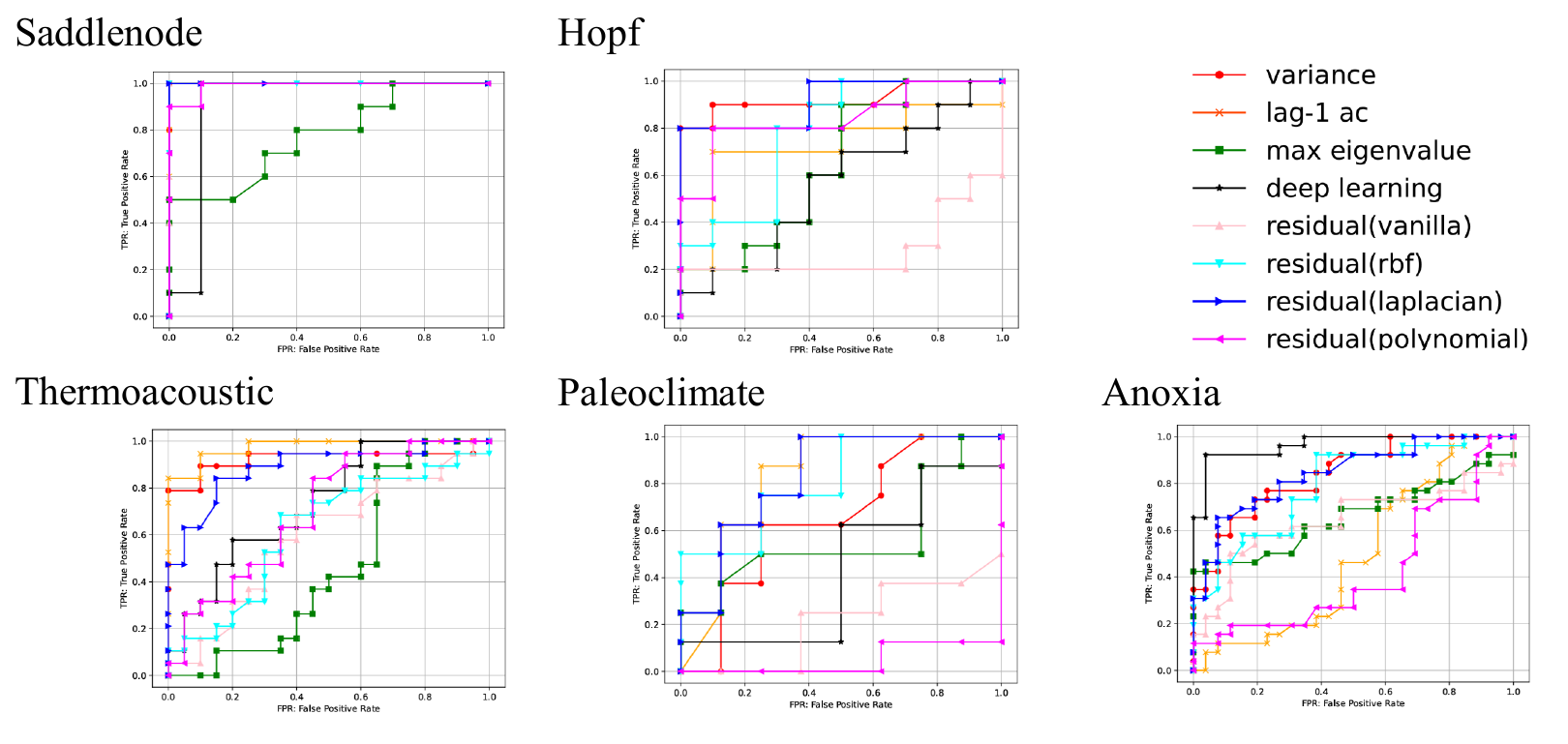}
    \caption{ROC curve (The Results for Saddle-node, Hopf, Thermoacoustic, Paleoclimate and Anoxic datasets are shown from top left to right.)}
    \label{fig:roc_all}
\end{figure*}

%\begin{table*}[h]
%    \centering
%    \caption{\red{Average of AUC}}
%    \label{tab:ave_auc}
%    \begin{tabular}{|c|c|c|c|}
%        \hline \hline
%         variance & lag-1 ac & max eigenvalue & deep learning \\ \hline
%         0.8700042822 & 0.805260939 & 0.6239425997 & 0.7215373715 \\ \hline \hline
%         residual(vanilla) & residual(rbf) & residual(laplacian) & residual(polynomial) \\ \hline
%         0.546326709 & 0.8019646528 & 0.8963067191 & 0.6016319293 \\ \hline \hline
%    \end{tabular}
%\end{table*}

We evaluate the proposed EWS described in Sections 3 and 4. The most important aspect of EWSs performance is achieving both a high detection rate and a low false-positive rate. Therefore, to compare the performance between the proposed and conventional EWSs, we prepare two artificial and three empirical data both for the tipping/non-tipping cases, and construct a receiver operating characteristic (ROC) curve by plotting the true positive rate (TPR) on the vertical axis against the false positive rate (FPR) on the horizontal axis. As conventional EWSs, we select variance\cite{ives1995measuring}, lag-1 autocorrelation\cite{dakos2008slowing,held2004detection}, maximum eigenvalue of the DMD matrix\cite{DONOVAN2022127152}, and deep learning method\cite{bury2021deep}. As an indicator for determining whether tipping occurs or not, we employ Kendall’s rank correlation coefficient between time and EWSs. This metric is appropriate when the bifurcation parameter is assumed to vary monotonically with time. (It should be noted that the parameters do not necessarily change monotonically in real-world situations.) The window size is set to half the total number of data points. For the DMD computation, the data is first embedded in a $400$-dimensional space using the delay coordinate. Also, when computing ResKMD with kernel EDMD, we employ three kernel functions: the radial basis function (RBF) kernel $k_{\mathrm{rbf}}(\mathbf{x}, \mathbf{x}') = \mathrm{exp}(-\gamma \|\mathbf{x}-\mathbf{x}'\|_2^2)$, the Laplacian kernel $k_{\mathrm{Laplacian}}(\mathbf{x}, \mathbf{x}') = \mathrm{exp}(-\gamma \|\mathbf{x}-\mathbf{x}'\|_1^2)$, and the polynomial kernel $k_{\mathrm{poly}}(\mathbf{x}, \mathbf{x}') = (\gamma \langle \mathbf{x}, \mathbf{x}' \rangle_2 + 1.0)^d$ with appropriate hyperparameter. The hyperparameter of the RBF and the Laplacian kernel is selected from $\gamma = 0.1$, $0.01$, $0.005$, $0.001$, $0.0005$, $0.0001$ by minimizing the ResKMD at the first time window. Also, The hyperparameter of the polynomial kernel is selected from $\gamma = 1.0$, $0.1$, $0.01$ and $d = 2.0$, $3.0$, $4.0$.

\begin{itemize}
\item Saddle-node and Subcritical Hopf Bifurcation

The local bifurcation introduced in Section 2 is known to retain their essential characteristics in a one- or two-dimensional dynamical system through the center manifold Theorem. Now, we focus on two local bifurcation, saddle-node bifurcation and subcritical Hopf bifurcation\cite{kuznetsov1998elements}. The normal form of saddle-node bifurcation is given by
\begin{align}
    \dot{x}=-\beta-x^2. \label{eq:normal_saddlenode}
\end{align}
When $\beta < 0$, the system has a stable fixed point at $x = \sqrt{-\beta}$ and an unstable fixed point at $x = -\sqrt{-\beta}$. As $\beta$ approaches $0$, these two fixed points move closer together and collide at $\beta = 0$, after which the system becomes unstable. Consequently, if the bifurcation parameter $\beta$ gradually increases over time, one observes how the equilibrium slowly shifts toward $0$ and eventually loses stability. Also, the normal form of subcritical Hopf bifurcation is given by
\begin{equation}
    \begin{aligned}
        \dot{x} &= \beta x - y - x(x^2 + y^2)(x^2 + y^2 - 1), \\
        \dot{y} &= x + \beta y - y(x^2 + y^2)(x^2 + y^2 -1).
    \end{aligned}
    \label{eq:normal_hopf}
\end{equation}
When $\beta < 0$, the system has a stable fixed point at the origin and an unstable limit cycle near the origin. On the other hand, $\beta > 0$, the origin becomes unstable, and a global stable limit cycle appears.

To examine the performance of our EWS, we generate two artificial data related to the saddle-node and subcritical Hopf bifurcation and calculate the ResKMD by applying Algorithms \ref{alg:ews_reskmd1}, \ref{alg:ews_reskmd2}. Time-series data of saddle-node bifurcation is generated by
\begin{align}
    \dot{x} = -(x+1)\big((x-1)^2-\beta\big), \label{eq:eq_saddlenode}
\end{align}
where the initial state $x(0) = 1.8$ and initial bifurcation parameter $\beta = 1.0$. We prepare some time-series in which the parameter gradually decreases at different rate. Time-series data of subcritical Hopf bifurcation is generated by
\begin{equation}
    \begin{aligned}
        \dot{x} &= \beta x - 2\pi y - x(x^2 + y^2)(x^2 + y^2 - 1), \\
        \dot{y} &= 2\pi x + \beta y - y(x^2 + y^2)(x^2 + y^2 -1).
    \end{aligned}
    \label{eq:eq_hopf}
\end{equation}
where the initial state $\big(x(0), y(0)\big) = (0.1, 0.0)$ and initial bifurcation parameter $\beta = -1.0$. Also, we prepare some time-series in which the parameter gradually increases at different rate. (We show an example in Fig. \ref{fig:ex_hopf}.) These time-series data are given by the explicit Runge-Kutta method of order $5$.

\item Tipping for Thermoacoustic System

Pavithran et al.\cite{pavi2021effect} investigated the subcritical Hopf bifurcation in a thermoacoustic system and discussed how the performance of EWSs varies depending on the rate at which the bifurcation parameter is changed.

\item Prediction of Climate Change

In Dakos' study\cite{dakos2008slowing}, it was reported that some EWSs can be detected for paleoclimate transitions, such as desertification in North Africa and warming events during the ice ages. In this experiment,  because of the limited availability of data, we increase the number of data points by means of spline interpolation. In principle, interpolation should be avoided in the computation of EWSs because it may introduce spurious correlations.

\item Marine Anoxic Events

Due to environmental factors such as global warming and eutrophication, the marine system in the eastern Mediterranean can sometimes abruptly transition into a state of oxygen depletion. Hennekam et al.\cite{hennekam2020early} have shown that some EWSs of this abrupt shift can be identified from sediments present in the oceanic region. In this experiment, we use the concentration of molybdenum (Mo) and uranium (U) in the sediments to calculate EWS, spline interpolation to increase the number of data points.
\end{itemize}

The results are shown in Fig. \ref{fig:roc_all}. Specifically, against saddle-node bifurcation, various EWSs achieve a clear separation between positive (showing tipping) and negative (not showing tipping) examples. Detecting a subcritical Hopf bifurcation is more challenging because these systems exhibit a pronounced delayed bifurcation. However, it can be seen that variance and ResKMD with the Laplacian kernel are good indicators for detecting tipping in this dataset. The remaining results show that the best performing EWS strongly depends on the dataset: lag-1 autocorrelation, ResKMD with the Laplacian kernel, and the deep learning method attain the highest detection rates for the thermoacoustic, paleoclimate, and anoxic datasets, respectively. From another perspective, across all datasets, ResKMD with the RBF and Laplacian kernel delivers consistently robust detection. Whereas the deep learning method relies on extensive pre-training to learn the systems near the tipping point, ResKMD with the RBF and Laplacian kernel reach comparable performance without any pre-training by accurately representing the latent systems close to the bifurcation point. Variance, lag-1 autocorrelation, and maximum eigenvalue of the DMD matrix likewise do not require pre-training.

\begin{figure*}[t]
    \centering
    \includegraphics[width=\textwidth]{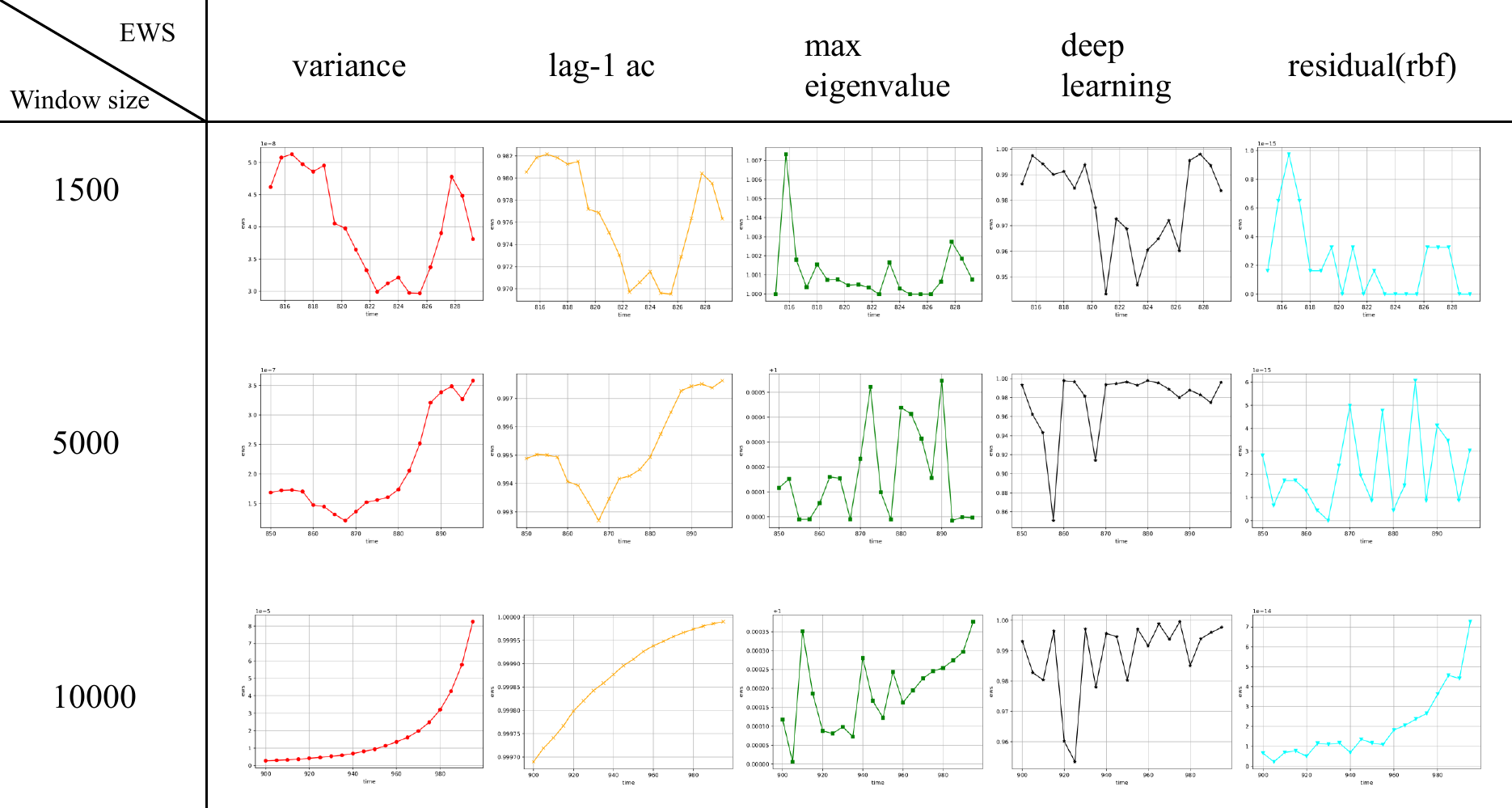}
    \caption{Behavior of EWSs against $1500$, $5000$, $10000$ Window Size\: This figure represents experimental results based on simulation data exhibiting saddle–node bifurcation.}
    \label{fig:vary_window}
\end{figure*}

\begin{figure*}[t]
    \centering
    \includegraphics[width=\textwidth]{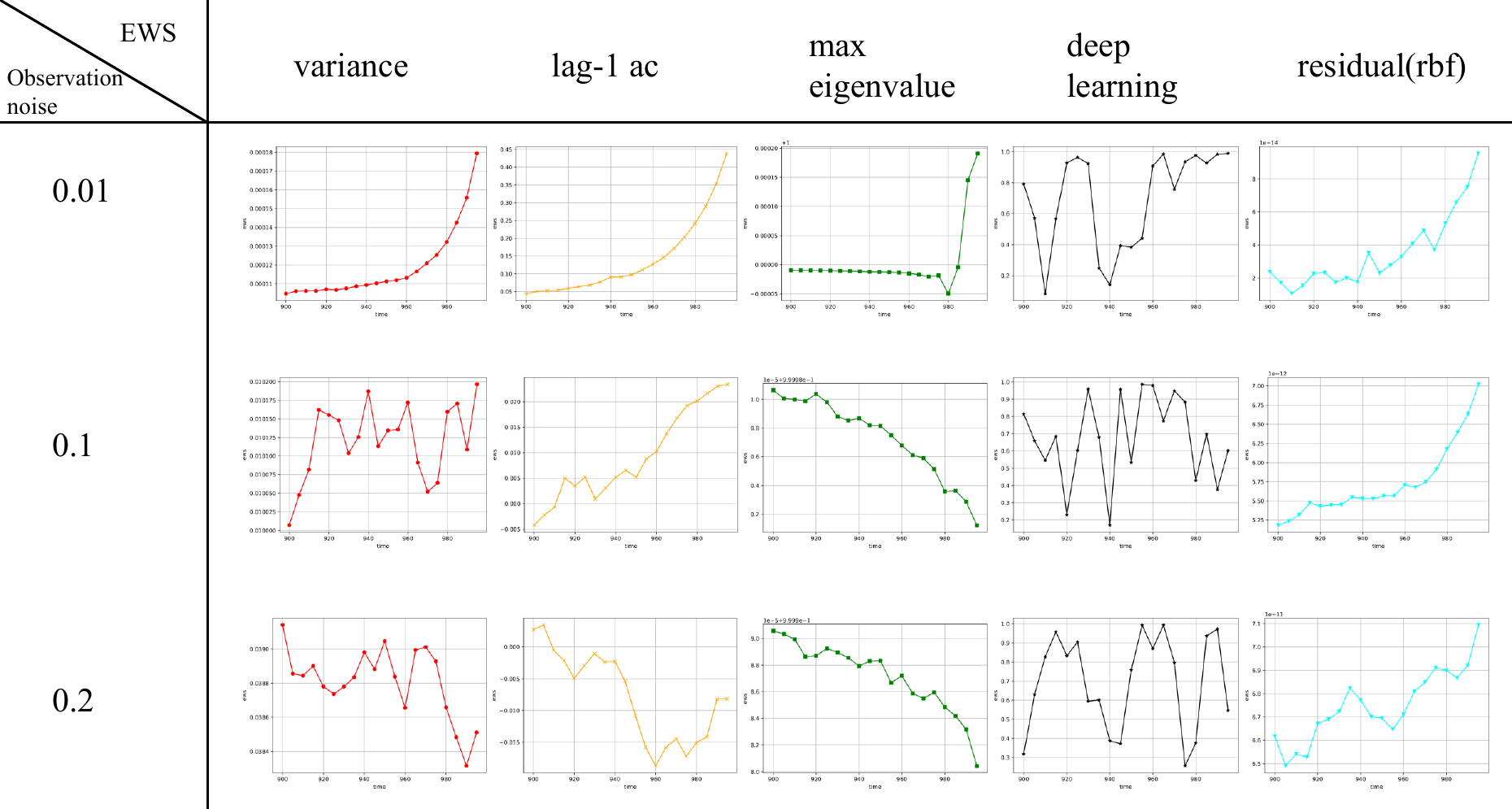}
    \caption{Behavior of against Observation Noise $\sigma_{\mathrm{obs}}^2 = 0.01, 0.1, 0.2$: This figure also represents experimental results based on simulation data exhibiting saddle–node bifurcation.}
    \label{fig:vary_obsnoise}
\end{figure*}

Other experiment results are shown in Fig. \ref{fig:vary_window} and Fig. \ref{fig:vary_obsnoise}. Fig. \ref{fig:vary_window} illustrates the variation of ResKMD with the RBF kernel when the window size is set to $1500$, $5000$, and $10000$ against the artificial data related to saddle-node bifurcation . In this experiment, the data are embedded into an $800$‐dimensional space using delay coordinates, and the kernel hyperparameter is tuned according to the procedure described earlier. What this result indicates is that ResKMD requires a large amount of data for computation. This is consistent with the well-known result experimentally demonstrated by Korda et al.\cite{korda2018convergence}, namely that a sufficiently large number of snapshot data is required in DMD to accurately approximate the spectra of the Koopman operator. However, this limitation of ResKMD is inconvenient in the context of early warning of tipping phenomena. For instance, in the paleoclimate or anoxia data, the available data may not be sufficient to detect the tendencies of tipping phenomena. Although we deal with this problem via spline interpolation in the experiment shown in Fig. \ref{fig:roc_all}, more carefully considered solutions are still required. Moreover, in kernel EDMD, it should be noted that increasing the window size leads to the “curse of dimensionality”. 

The experiment presented in Fig. \ref{fig:vary_obsnoise} investigates how performance changes when observation noise is added to simulation data exhibiting saddle–node bifurcation. The observation noise is generated based on i.i.d. Gaussian white noise with density $\mathcal{N}(0,\sigma^2_{\mathrm{obs}})$ where $\sigma_{\mathrm{obs}} > 0$. This result confirms that ResKMD with the RBF kernel serves as a robust EWS against observation noise. On the other hand, variance and lag-1 autocorrelation's performance deteriorates significantly when observation noise is strongly mixed, because variance and lag-1 autocorrelation directly measure stochastic resilience in the data space. The deep learning method also do not show good performance against data with strong observation noise. This method can be interpreted as outputting the probability that the observed behavior resembles that of tipping phenomena present in the train data. However, since the data are contaminated by observation noise (strictly speaking, exhibiting behavior not included in the train data), accurate prediction becomes difficult. This is a weakness inherent to pre‐training approaches.

\section{Conclusions}
In this paper, we newly propose an EWS based on the Koopman operator. Our EWS is developed by concerning the spectral property of the Koopman operator near the bifurcation point. As the absolute value of the dominant Koopman eigenvalue approaches $1$, the estimated error between the true Koopman operator and the proxy of the Koopman operator by its mode decomposition increases because of its continuous spectrum. This mathematical fact shows that ResKMD behaves as an EWS. In particular, we have proved that it can be considered as generalized stochastic resilience and strictly increases as the bifurcation point approaches. Moreover, it can be numerically computed from data via ResDMD, and we demonstrate its effectiveness through experiments on five synthetic and real datasets. As noted at the end of Section 5, the proposed method exhibits robust detection performance across diverse data types.

Although our results are promising, several extensions remain. First, we can change the method to estimate the Koopman operator. For example, neural DMD\cite{iwata2020neural} and jet DMD\cite{ishikawa2024koopman} are able to estimate the Koopman operator more accurately than exact DMD and kernel EDMD. By employing these methods, signs of the tipping point are expected to be detected with greater accuracy at an earlier stage. Second, it is still an open challenge to devise a versatile EWS capable of detecting global bifurcation-induced tipping\cite{kuznetsov1998elements}, such as transitions to chaotic attractors, as well as rate-induced tipping\cite{wieczorek2023rate}, which depends on the rate of parameter variation.

\begin{acknowledgements}
This work was partially supported by JSPS KAKENHI Grant Numbers JP22H00516, JP22H05106 and JST CREST Grant Number JPMJCR1913. In addition, the authors thank the anonymous reviewers and the editor for their constructive comments and valuable suggestions, which greatly helped to improve the quality and clarity of this article.
\end{acknowledgements}

{\small
\noindent\textbf{Data Availability} The three empirical datasets \cite{pavi2021effect,dakos2008slowing,hennekam2020early} and the EWS method proposed by Bury et al. \cite{bury2021deep} were obtained and, in part, directly used as provided in \url{https://github.com/ThomasMBury/deep-early-warnings-pnas?tab=readme-ov-file}. All codes and data used in this study are available in this GitHub repository \url{https://github.com/uosaka-mlsyslab/reskmd_for_ews_2025}.

\noindent{}A preprint of this work can be found at~\cite{2508.1965}.

\noindent\textbf{Conflict of Interest} The authors declare that they have no conflict of interest.}

\appendix
\section{Proof of Corollary \ref{cor01}}
For a fixed $\mathbf{x} \in X$, the expectation of $\bigl\| \mathbf{g}_1(\mathbf{F}_{\boldsymbol{\omega}}(\mathbf{x}; \beta)) + \mathbf{g}_2( \mathbf{x})\bigr\|_2^2$ is expanded as follows, 
\begin{align}
    &\mathrm{E}\Bigl[ \bigl\| \mathbf{g}_1\circ\mathbf{F}_{\boldsymbol{\omega}}(\mathbf{x}; \beta) + \mathbf{g}_2( \mathbf{x})\bigr\|_2^2 \Bigr] \notag \\
    & \qquad = \mathrm{E}\Bigl[ \bigl\| \mathbf{g}_1\circ\mathbf{F}_{\boldsymbol{\omega}}(\mathbf{x}; \beta)\bigr\|_2^2 \Bigr] + 2\mathrm{Re}\big\langle U_{(1)}\mathbf{g}_1(\mathbf{x}), \mathbf{g}_2(\mathbf{x}) \big\rangle_2 + \|\mathbf{g}_2(\mathbf{x})\|_2^2 \notag \\
    & \qquad =\bigl\| U_{(1)}\mathbf{g}_1(\mathbf{x}) + \mathbf{g}_2(\mathbf{x}) \bigr\|_2^2 \notag \\
    & \qquad \qquad+ \mathrm{E}\Bigl[ \bigl\| \mathbf{g}_1\circ\mathbf{F}_{\boldsymbol{\omega}}(\mathbf{x}; \beta) \bigr\|_2^2 \Bigr] - \Bigl\| \mathrm{E}\bigl[ \mathbf{g}_1\circ\mathbf{F}_{\boldsymbol{\omega}}(\mathbf{x}; \beta) \bigr] \Bigr\|_2^2 \notag \\
    & \qquad =\bigl\| U_{(1)}\mathbf{g}_1(\mathbf{x}) + \mathbf{g}_2(\mathbf{x}) \bigr\|_2^2 + \operatorname{Tr} \bigg( \mathrm{Cov}[\mathbf{g}_1 \circ \mathbf{F}_{\boldsymbol{\omega}}(\mathbf{x}; \beta)] \bigg). \notag
\end{align}
By integrating over $\mathbf{x}$ with respect to the measure $\mu$, we can obtain the result \eqref{eq:multi_biasvariance}.

\section{Proof of Proposition \ref{prop01}}
We assume that the eigenfunctions form an orthonormal basis in $L^2(X, \mu)$, we have
\begin{align}
    \bigl\|U_{(1)}\mathbf{g}-U_{(1),\Phi_m}\mathbf{g}\bigr\|_{L^2(X, \mu)^M}^2 = \sum_{n = m+1}^\infty \lvert \lambda_n \rvert^2 \| \mathbf{v}_{\lambda_n} \|_2^2. \label{A:residual_kmd}
\end{align}
Here, the modes are given by $\mathbf{v}_{\lambda_n} = \langle \mathbf{g},\phi_{\lambda_n} \rangle$ and are satisfied 
\begin{align}
    \| \mathbf{g} \|_{L^2(X, \mu)^M}^2 = \sum_{n=1}^{\infty} \lvert \langle \mathbf{g}, \phi_{\lambda_n} \rangle \rvert^2. \label{A:obs_norm}
\end{align}
Therefore, we can bound ResKMD from the above as
\begin{align}
    &\bigl\|U_{(1)}\mathbf{g}-U_{(1),\Phi_m}\mathbf{g}\bigr\|_{L^2(X, \mu)^M}^2 \notag \\
    & \qquad\le \lvert \lambda_{m+1} \rvert^2 \Bigl( \|\mathbf{g}\|_{L^2(X, \mu)^M}^2 - \sum_{n=1}^{m} \lvert \langle \mathbf{g}, \phi_{\lambda_n} \rangle \rvert^2 \Bigr) \notag \\
    & \qquad \le \lvert \lambda_{m+1} \rvert^2 \|\mathbf{g}\|_{L^2(X, \mu)^M}^2. \label{A:reskmd_upper}
\end{align}
In addition, the lower bound can be established:
\begin{align}
    \lvert \lambda_{m+1} \rvert^2 \| \mathbf{v}_{\lambda_{m+1}} \|_2^2 &\le \bigl\|U_{(1)}\mathbf{g}-U_{(1),\Phi_m}\mathbf{g}\bigr\|_{L^2(X, \mu)^M}^2. \label{A:reskmd_lower}
\end{align}

\section{The Proof of Eq. \eqref{eq:reskmd_aprox}}
We present here the derivation of Eq. \eqref{eq:reskmd_aprox}. By definition, we have
\begin{align}
    &U_{(1), \Phi_m} \mathbf{g}(\mathbf{x}) = \sum_{n=1}^m \lambda_{n} \mathbf{v}_{\lambda_n} \phi_{\lambda_n}(\mathbf{x}) \label{eqB101} \\
    &\mathbf{g}(\mathbf{x}) = \sum_{n=1}^m \mathbf{v}_{\lambda_n} \phi_{\lambda_n}(\mathbf{x}). \label{B:restrictedkoop_observable}
\end{align}
So, the error between the true Koopman operator $U_{(1)}$ and the restricted Koopman operator $U_{(1), \Phi_m}$ is
\begin{align}
    U_{(1)}\mathbf{g} - U_{(1), \Phi_m}\mathbf{g} = \sum_{n=1}^m \mathbf{v}_{\lambda_n} (U_{(1)}\phi_{\lambda_n}-\lambda_n \phi_{\lambda_n}). \label{B:resid_byeig}
\end{align}
Hence, considering the norm of error \eqref{B:resid_byeig} under $\mu$, we may exchange the integral and the summation if $U_{(1)}\phi_{\lambda_n} - \lambda_n\phi_{\lambda_n}$ are mutually orthogonal, leading to
\begin{align}
    &\int_X |U_{(1)}\mathbf{g} - U_{(1), \Phi_m}\mathbf{g}|^2 d\mu(\mathbf{x}) \notag \\
    & \qquad =\sum_{n=1}^m |\mathbf{v}_{\lambda_n}|^2 \int_X | U_{(1)}\phi_{\lambda_n}(\mathbf{x})-\lambda_n \phi_{\lambda_n}(\mathbf{x}) |^2 d\mu(\mathbf{x}) \label{B:residnorm_another}
\end{align}
Therefore, ResKMD is computed as
\begin{align}
    \mathrm{res}[U_{(1)}, U_{(1), \Phi_m}]^2 &= \sum_{n=1}^m |\mathbf{v}_{\lambda_n}|^2 \int_X | U_{(1)}\phi_{\lambda_n}(\mathbf{x})-\lambda_n \phi_{\lambda_n}(\mathbf{x}) |^2 d\mu(\mathbf{x}) \notag \\
    &= \sum_{n=1}^m |\mathbf{v}_{\lambda_n}|^2 \mathrm{res}[\lambda_n, \phi_{\lambda_n}]^2, \label{B: reskmd_by_resdmd}
\end{align}
where in the last step, we have assumed that each $\phi_{\lambda_n}$ is normalized. Proceeding analogously to Proposition \ref{prop01}, we observe that the norm of the mode remains nearly constant. Hence, Eq. \eqref{eq:reskmd_aprox} follows.

\bibliographystyle{spmpsci}      % mathematics and physical sciences
\bibliography{reference}   % name your BibTeX data base

\begin{thebibliography}{10}
\providecommand{\url}[1]{{#1}}
\providecommand{\urlprefix}{URL }
\expandafter\ifx\csname urlstyle\endcsname\relax
  \providecommand{\doi}[1]{DOI~\discretionary{}{}{}#1}\else
  \providecommand{\doi}{DOI~\discretionary{}{}{}\begingroup \urlstyle{rm}\Url}\fi

\bibitem{alfatlawi2019incremental}
Alfatlawi, M., Srivastava, V.: An incremental approach to online dynamic mode decomposition for time-varying systems with applications to eeg data modeling.
\newblock arXiv preprint arXiv:1908.01047  (2019)

\bibitem{brunton2022modern}
Brunton, S.L., Budi\v{s}i\'{c}, M., Kaiser, E., Kutz, J.N.: Modern koopman theory for dynamical systems.
\newblock SIAM Review \textbf{64}(2), 229--340 (2022)

\bibitem{bury2021deep}
Bury, T.M., Sujith, R.I., Pavithran, I., Scheffer, M., Lenton, T.M., Anand, M., Bauch, C.T.: Deep learning for early warning signals of tipping points.
\newblock Proceedings of the National Academy of Sciences \textbf{118}(39), e2106140118 (2021)

\bibitem{colbrook2023residual}
Colbrook, M.J., Ayton, L.J., Sz{\H{o}}ke, M.: Residual dynamic mode decomposition: robust and verified koopmanism.
\newblock Journal of Fluid Mechanics \textbf{955}, A21 (2023)

\bibitem{colbrook2024beyond}
Colbrook, M.J., Li, Q., Raut, R.V., Townsend, A.: Beyond expectations: residual dynamic mode decomposition and variance for stochastic dynamical systems.
\newblock Nonlinear Dynamics \textbf{112}(3), 2037--2061 (2024)

\bibitem{vcrnjaric2020koopman}
{\v{C}}rnjari{\'c}-{\v{Z}}ic, N., Ma{\'c}e{\v{s}}i{\'c}, S., Mezi{\'c}, I.: Koopman operator spectrum for random dynamical systems.
\newblock Journal of Nonlinear Science \textbf{30}, 2007--2056 (2020)

\bibitem{dakos2008slowing}
Dakos, V., Scheffer, M., Van~Nes, E.H., Brovkin, V., Petoukhov, V., Held, H.: Slowing down as an early warning signal for abrupt climate change.
\newblock Proceedings of the National Academy of Sciences \textbf{105}(38), 14308--14312 (2008)

\bibitem{DONOVAN2022127152}
Donovan, G., Brand, C.: Spatial early warning signals for tipping points using dynamic mode decomposition.
\newblock Physica A: Statistical Mechanics and its Applications \textbf{596}, 127152 (2022).
\newblock \doi{https://doi.org/10.1016/j.physa.2022.127152}

\bibitem{gaspard1995spectral}
Gaspard, P., Nicolis, G., Provata, A., Tasaki, S.: Spectral signature of the pitchfork bifurcation: Liouville equation approach.
\newblock Physical Review E \textbf{51}(1), 74 (1995)

\bibitem{held2004detection}
Held, H., Kleinen, T.: Detection of climate system bifurcations by degenerate fingerprinting.
\newblock Geophysical Research Letters \textbf{31}(23) (2004)

\bibitem{hennekam2020early}
Hennekam, R., van~der Bolt, B., van Nes, E.H., de~Lange, G.J., Scheffer, M., Reichart, G.J.: Early-warning signals for marine anoxic events.
\newblock Geophysical Research Letters \textbf{47}(20), e2020GL089183 (2020)

\bibitem{ishikawa2024koopman}
Ishikawa, I., Hashimoto, Y., Ikeda, M., Kawahara, Y.: Koopman operators with intrinsic observables in rigged reproducing kernel hilbert spaces.
\newblock arXiv preprint arXiv:2403.02524  (2024)

\bibitem{ives1995measuring}
Ives, A.R.: Measuring resilience in stochastic systems.
\newblock Ecological Monographs \textbf{65}(2), 217--233 (1995)

\bibitem{iwata2020neural}
Iwata, T., Kawahara, Y.: Neural dynamic mode decomposition for end-to-end modeling of nonlinear dynamics.
\newblock arXiv preprint arXiv:2012.06191  (2020)

\bibitem{koopman1931hamiltonian}
Koopman, B.O.: Hamiltonian systems and transformation in hilbert space.
\newblock Proceedings of the National Academy of Sciences \textbf{17}(5), 315--318 (1931)

\bibitem{korda2018convergence}
Korda, M., Mezi{\'c}, I.: On convergence of extended dynamic mode decomposition to the koopman operator.
\newblock Journal of Nonlinear Science \textbf{28}(2), 687--710 (2018)

\bibitem{kuznetsov1998elements}
Kuznetsov, Y.A., Kuznetsov, I.A., Kuznetsov, Y.: Elements of applied bifurcation theory, vol. 112.
\newblock Springer (1998)

\bibitem{lenton2008tipping}
Lenton, T.M., Held, H., Kriegler, E., Hall, J.W., Lucht, W., Rahmstorf, S., Schellnhuber, H.J.: Tipping elements in the earth's climate system.
\newblock Proceedings of the national Academy of Sciences \textbf{105}(6), 1786--1793 (2008)

\bibitem{mauroy2016global}
Mauroy, A., Mezi{\'c}, I.: Global stability analysis using the eigenfunctions of the koopman operator.
\newblock IEEE Transactions on Automatic Control \textbf{61}(11), 3356--3369 (2016)

\bibitem{mauroy2020koopman}
Mauroy, A., Susuki, Y., Mezic, I.: Koopman operator in systems and control, vol. 484.
\newblock Springer (2020)

\bibitem{mezic2005spectral}
Mezi{\'c}, I.: Spectral properties of dynamical systems, model reduction and decompositions.
\newblock Nonlinear Dynamics \textbf{41}(1), 309--325 (2005)

\bibitem{2508.1965}
Miyauchi, Y., Ikeda, M., Kawahara, Y.: Generalized stochastic resilience for early warning signals based on koopman operator (2025)

\bibitem{pavi2021effect}
Pavithran, I., Sujith, R.I.: Effect of rate of change of parameter on early warning signals for critical transitions.
\newblock Chaos: An Interdisciplinary Journal of Nonlinear Science \textbf{31}(1), 013116 (2021).
\newblock \doi{10.1063/5.0025533}

\bibitem{rockstrom2009safe}
Rockstr{\"o}m, J., Steffen, W., Noone, K., Persson, {\AA}., Chapin, F.S., Lambin, E.F., Lenton, T.M., Scheffer, M., Folke, C., Schellnhuber, H.J., et~al.: A safe operating space for humanity.
\newblock nature \textbf{461}(7263), 472--475 (2009)

\bibitem{scheffer2009early}
Scheffer, M., Bascompte, J., Brock, W.A., Brovkin, V., Carpenter, S.R., Dakos, V., Held, H., Van~Nes, E.H., Rietkerk, M., Sugihara, G.: Early-warning signals for critical transitions.
\newblock Nature \textbf{461}(7260), 53--59 (2009)

\bibitem{takeishi2017subspace}
Takeishi, N., Kawahara, Y., Yairi, T.: Subspace dynamic mode decomposition for stochastic koopman analysis.
\newblock Physical Review E \textbf{96}(3), 033310 (2017)

\bibitem{takens2006detecting}
Takens, F.: Detecting strange attractors in turbulence.
\newblock In: Dynamical Systems and Turbulence, Warwick 1980: proceedings of a symposium held at the University of Warwick 1979/80, pp. 366--381. Springer (2006)

\bibitem{van2007slow}
Van~Nes, E.H., Scheffer, M.: Slow recovery from perturbations as a generic indicator of a nearby catastrophic shift.
\newblock The American Naturalist \textbf{169}(6), 738--747 (2007)

\bibitem{wieczorek2023rate}
Wieczorek, S., Xie, C., Ashwin, P.: Rate-induced tipping: Thresholds, edge states and connecting orbits.
\newblock Nonlinearity \textbf{36}(6), 3238 (2023)

\bibitem{williams2015data}
Williams, M.O., Kevrekidis, I.G., Rowley, C.W.: A data--driven approximation of the koopman operator: Extending dynamic mode decomposition.
\newblock Journal of Nonlinear Science \textbf{25}, 1307--1346 (2015)

\bibitem{wissel1984universal}
Wissel, C.: A universal law of the characteristic return time near thresholds.
\newblock Oecologia \textbf{65}, 101--107 (1984)

\bibitem{zhang2019online}
Zhang, H., Rowley, C.W., Deem, E.A., Cattafesta, L.N.: Online dynamic mode decomposition for time-varying systems.
\newblock SIAM Journal on Applied Dynamical Systems \textbf{18}(3), 1586--1609 (2019)

\end{thebibliography}

\end{document}